\documentclass[11pt,a4paper,twoside]{article}
\usepackage{etoolbox}  % used in jmpag.sty
\usepackage{forloop}    % used in jmpag.sty
\usepackage{amsthm,amsfonts,amsmath,amscd,amssymb}
\usepackage{latexsym}
\usepackage{euscript}
\usepackage{enumitem}
\usepackage{graphicx}
\usepackage{tikz}
\usepackage{caption}
\usepackage{subcaption}
\usepackage{cite}
\usepackage[utf8]{inputenc}
\usepackage[english]{babel}
\usepackage{xcolor}
\definecolor{darkgreen}{HTML}{00A000}
\definecolor{darkblue}{HTML}{0000A0}
\definecolor{darkred}{HTML}{D00000}
\usepackage[colorlinks=true]{hyperref}%%%%%%%%%[,backref]
\hypersetup{linkcolor=darkred, urlcolor=darkblue, citecolor=darkgreen}%%%%%%%%
\usepackage{jmpag}

\begin{document}

%%%%%%%%%%%%%% TITLE and AUTHORS %%%%%%%%%%%%%%%%%%

% The full title of your paper
\title[Controllability Problems for the Heat Equation with Variable Coefficients]%
{Controllability Problems for the Heat Equation with Variable Coefficients on a Half-Axis Controlled by the Neumann Boundary Condition}

% If the title of your paper is too long, use the short title please. It should not exceed 60 symbols (including blank spaces).

%\title[Short Title]{Title}

% Author1
%FirstName1 LastName1
\author{Larissa  Fardigola}
%Institution1, address1, City1, Postal Code1, Country
\address{B. Verkin Institute for Low Temperature Physics and Engineering of the
National Academy of Sciences of Ukraine, 47 Nauky Ave., Kharkiv, 61103, Ukraine,\\
V.N. Karazin Kharkiv National University, 4 Svobody Sq., Kharkiv, 61022}
%e-mail1
\email{fardigola@ilt.kharkov.ua}

% Author2
%FirstName2 LastName2
\author{Kateryna  Khalina}
%Institution2, address2, City2, Postal Code2, Country2
\address{B. Verkin Institute for Low Temperature Physics and Engineering of the
National Academy of Sciences of Ukraine, 47 Nauky Ave., Kharkiv, 61103, Ukraine}
%e-mail2
\email{khalina@ilt.kharkov.ua}

%%%%%% Received December 18, 2021 

%%%%%%%%%%%%%%%% END TITLE and AUTHORS %%%%%%%%%%%%%%%%%%%

\BeginPaper %%%%%%% do not remove this command

%%%%%%%%%%%%%%%%%%%%%% MY MACROS %%%%%%%%%%%%%%%%%%%%

\renewcommand\textfraction{0}

\renewcommand\topfraction{1}

\renewcommand\bottomfraction{1}

\newcommand{\LL}{L_{-1,(3)}^2(\mathbb R_+)}

%%%%%%%%%%%%%%%%%%%%%%%%%%%%%%%%%%%%%%%%%%%%%%%%%%%%%%%%%%%%%%%%%%%%%%%

%%%%%%%%%%%%%%%%%%%%%%
\newcommand{\fnl}{\left[\kern-0.1em\left]}
\newcommand{\fnr}{\right[\kern-0.1em\right]}
\newcommand{\al}{\left\langle\kern-0.3ex\left\langle}
\newcommand{\ar}{\right\rangle\kern-0.3ex\right\rangle}
\newcommand{\sgn}{\mathop{\mathrm{sgn}}}
\newcommand{\supp}{\mathop{\mathrm{supp}}}
%\newcommand{\diag}{\mathop{\mathrm{diag}}

%%%%%%%%%%%%%%%%%%%%%%%%%%%%%%%%%%%%%%%%%%%%%%%%%%%%%%%%%%%%%%%
%\newcommand{\zco}{{\mbox{\scriptsize$\bigcirc$\kern-6.1pt\raisebox{-0.5pt}{$0$}}}}
\newcommand{\zco}{{\text{$\bigcirc$\kern-6.2pt\raisebox{-0.5pt}{$0$}}}}
\newcommand{\zcx}[1]{~\mbox{#1\kern-7.5pt\raisebox{0.5pt}{$\bigcirc$}}}
%%%%%%%%%%%%%%%%%%%%%%%%%%%%%%%%%%%%%%%%%%%%%%%%%%%%%%%%%%%%%%%%%%%%%%%%%

%%%%%%%%%%%%%%%%%%%%% END MY MACROS %%%%%%%%%%%%%%%%%%%%%%

%%%%%%%%%%%%%%%%%%%%%%% ABSTRACT %%%%%%%%%%%%%%%%%%%%%%%%

\begin{abstract}

In the paper, the problems of controllability and approximate controllability are studied for the control system
$w_t=\frac{1}{\rho}\left(kw_x\right)_x+\gamma w$, $\left.\left(\sqrt{\frac{k}{\rho}}w_x\right)\right|_{x=0}=u$, $x>0$, $t\in(0,T)$, where $u$ is a control, $u\in L^\infty(0,T)$. It is proved that each initial state of the control system is approximately 
controllable to any  target state in a given time $T>0$. To obtain this result, the transformation operator generated 
by the equation data $\rho$, $k$, $\gamma$ is applied. The results are illustrated by examples.

\key{heat equation, controllability, approximate controllability}

\end{abstract}

%%%%%%%%%%%%%%%%%%%% END ABSTRACT %%%%%%%%%%%%%%%%%%%%%%%%%%%%%

%The title of your section \href{mailto:#1}{{\mdseries\ttfamily #1}}}
\section{Introduction}\label{sect1}

Controllability problems for the heat equation with constant and variable coefficients were studied in  a number of papers (see, e.g., \cite{Bic,  CMV, CN, DE, FKh, FKh2, FatRus, FCZ, IY, MV, MZua1, SMEZ, MP, SS,  ZXY, EZua}).
However, there are only a few papers where these problems were investigated for \emph{the heat equation with constant coefficients on unbounded domains} (see, e.g., \cite{CMV, FKh, FKh2, MZua1, SMEZ}), and it seems these problems were investigated for \emph{the heat equation with variable coefficients on unbounded domains} only in \cite{FKh3}.

The paper deals with the controllability problems for the heat equation with variable coefficients on a half-axis controlled by the Neumann boundary condition. Consider the following control system:
\begin{align}
&w_t=\frac{1}{\rho}\left(kw_x\right)_x+\gamma w,&& x\in(0,+\infty),\ t\in(0,T),\label{eq}\\ 
&\left.\left(\sqrt{\frac{k}{\rho}}w_x\right)\right|_{x=0}=u,&& t\in(0,T),\label{bc}\\ 
&w(\cdot,0)=w^0,&&x\in(0,+\infty).\label{ic}
\end{align}
Here $T>0$ is a  constant; $\rho$, $k$, $\gamma$,  and $w^0$ are given functions; $u\in L^\infty(0,T)$ is a control. We  assume $\rho, k\in C^1[0,+\infty)$ are positive on $[0,+\infty)$, $(\rho k)\in C^2[0,+\infty)$, $(\rho k)'(0)=0$. Consider the even extensions of $\rho$, $k$, $\gamma$. Throughout the paper we will denote these extensions by the same symbols $\rho$, $k$, $\gamma$, respectively. Denote
\begin{equation}\label{sigma}
\sigma(x)=\int_0^x\sqrt{\rho(|\xi|)/k(|\xi|)}\,d\xi,\quad x\in\mathbb R.
\end{equation}
We assume
\begin{equation}\label{sigma-a}
\sigma(x)\to+\infty\quad\text{as }x\to+\infty.
\end{equation}
Put $Q_2(\rho,k)=\sqrt{ k/\rho}\big( Q_1(\rho,k) \big)'
+\big(Q_1(\rho,k) \big)^2$, $Q_1(\rho,k)=\sqrt{ k/\rho}(k\rho)'/(4k\rho)$.
We also assume
\begin{equation}\label{pot1}
Q_2(\rho,k)-\gamma\in L^\infty(0,+\infty)\bigcap C^1[0,+\infty)
\end{equation}
and
\begin{equation}\label{pot2}
\sqrt{\frac \rho k}\left( Q_2(\rho,k)-\gamma\right)\sigma\in L^1(0,+\infty).
\end{equation}
We consider control system \eqref{eq}--\eqref{ic} in modified Sobolev spaces (see Section \ref{spaces}).

In \cite{FKh3}, controllability problems for the heat equation with variable coefficients on a half-axis controlled by the Dirichlet boundary condition are studied. The general methods applied in the present paper are similar to those from paper \cite{FKh3}. But for the case of the Neumann boundary condition,  different spaces and operators are used that  caused  different technique of proofs of main results.

Theorems \ref{lackne} and \ref{thapprne} (see Section \ref{spaces} below) are the main result of the paper.  It is proved that each initial state of the control system is approximately controllable to any target state in a given time $T>0$ (Theorem \ref{thapprne}). In the case of constant coefficients $(\rho=k=1, \gamma=0)$, the result of this theorem has been obtained earlier in \cite{FKh2}. In the case of variable coefficients, this result is similar to those of papers \cite{F2, F1, disF, F3} for the wave equation with variable coefficients on a half-axis controlled either by the Dirichlet or by the Neumann boundary condition. However, the methods for obtaining the results are essentially different because of entirely different nature of the heat and wave equations. They are compared below. If an initial state of control system is controllable to the origin  then the initial state is also the origin (see Theorem \ref{lackne}). In the case of constant coefficients $(\rho=k=1, \gamma=0)$, the result of this theorem has been obtained earlier in \cite{FKh2}. This result is similar to that of the paper \cite{SMEZ}.

To study control system \eqref{eq}--\eqref{ic}, we use the transformation operator $\widehat{\mathbb T}$ and the modified Sobolev spaces $\widehat{\mathbb H}^s$, $s=\overline{-1,1}$. This operator $\widehat{\mathbb T}: \widehat H^{-1}\to \widehat{\mathbb H}^{-1}$ together with the spaces $\widehat{\mathbb H}^s$, $s=\overline{-1,1}$, associated with the equation data $(\rho,k,\gamma)$ are introduced and studied in \cite{F2, F1, disF, F3}. The definitions of $\widehat{\mathbb T}$, $\widehat{\mathbb H}^s$, and $\widehat H^s$ are given below in Section \ref{spaces}.

The operator $\widehat{\mathbb T}$ is a continuous one-to-one mapping between the spaces $\widehat H^s$ and $\widehat{\mathbb H}^s$. Moreover, it is one-to-one mapping between the set of the solutions to \eqref{eq}--\eqref{ic}  with constant coefficients ($\rho=k=1$, $\gamma=0$) where $u = u^{110}\in L^\infty(0,T)$ and the set of the solutions to this problem with variable coefficients $\rho, k, \gamma$ where $u=u^{\rho k\gamma}\in L^\infty(0,T)$ (see below Theorems \ref{thtransf1ne} and \ref{thtransf2ne}). Note that $u^{110}$ and $u^{\rho k\gamma}$ are different generally speaking.
The proofs of the main results of the paper are based on Theorems \ref{thtransf1ne} and \ref{thtransf2ne} proved in Section \ref{operators}. The control system with variable coefficients $\rho, k, \gamma$ replicates the controllability properties of the control system with constant coefficients ($\rho=k=1$, $\gamma=0$) and vice versa. 

The last result also holds true for the wave equation on a half-axis \cite{F2, F1, disF, F3}. But the proofs are essentially different for the cases of the wave and heat equations. Applying the operator $\widehat{\mathbb T}^{-1}$ to a solution to the equation with variable coefficients $\rho,k,\gamma$ and a control $u=u^{\rho k\gamma}\in L^\infty(0,T)$, we obtain a solution to the equation with the constant coefficients $\rho=k=1$, $\gamma=0$ and a control $u = u^{110}\in L^\infty(0,T)$ different from the control $u^{\rho k\gamma}$. To find and to estimate the control $u^{110}$, we have to solve an integral equation of the form
\begin{equation}
\label{inteq}
u^{110}(t)=f(t)+\int_0^t P(t-\xi) u^{110}(\xi) \, d\xi,\quad t\in[0,T].
\end{equation}
In the case of the wave equation, it has been proved that $f$ and $P$ are bounded on $[0,T]$ \cite{F2, F1, disF, F3}. Therefore, the integral operator in the right-hand side of \eqref{inteq} is of the Hilbert--Schmidt type. Hence, the Fredholm alternative together with the generalized Gronwall theorem can be applied  to solve \eqref{inteq} in $L^2(0,T)$ and estimate the solution $u^{110}$ in $L^\infty(0,T)$ when we deal with the wave equation \cite{F2, F1, disF, F3}. 
In the case of the heat equation, it has been proved that $f$ and $\sqrt{(\cdot)}P$ are bounded on $[0,T]$ (hence, $P(\xi)=O(1/\sqrt \xi)$ as $\xi\to0^+$) \cite{FKh3}. 
That is why the integral operator in the right-hand side of \eqref{inteq} is not of the Hilbert--Schmidt type, and the Fredholm alternative is not applicable  in the general case. 
The Banach fixed-point theorem is also not applicable in general case. That is why the method of successive approximations has been used  to construct a solution to \eqref{inteq} on $[0,T]$. Then  the Banach fixed-point theorem has been applied in $L^2$-space on small intervals to prove the uniqueness of the solution \cite{FKh3}. This result is recalled in  Lemma \ref{exist} below. 

Since the control system with variable coefficients $\rho, k, \gamma$ replicates the controllability properties of the control system with constant coefficients ($\rho=k=1$, $\gamma=0$), we obtain the controllability properties of the first control system from the controllability properties of the second one by applying the operator $\widehat{\mathbb T}$, i.e., we obtain Theorems \ref{lackne} and \ref{thapprne} by applying Theorems \ref{thtransf1ne} and \ref{thtransf2ne} in Section \ref{spaces}.

The obtained results are illustrated by examples in Section \ref{exxx}.

\section{Spaces, operators and main results}\label{spaces}

Let us give definitions of the spaces used in the paper.

Let $\Omega=(0,+\infty)$ or $\Omega=\mathbb R$.
Let   $\EuScript D(\Omega)$ be the space of finite infinitely differentiable functions whose support is finite and is contained in $\Omega$. 
For $\varphi\in L_{\text{loc}}^2(\Omega)$ we consider $\varphi'\in \EuScript D'(\Omega)$.

By $H^p$, $p=0,1$, denote  the Sobolev spaces
\begin{equation*}
H^p=\left\lbrace \varphi\in L^2_{\mathrm{loc}}(\mathbb{R})\mid \forall m=\overline{0,p}\ \varphi^{(m)}\in L^2(\mathbb{R})\right\rbrace 
\end{equation*}
with the norm
\begin{equation*}
\left\| \varphi\right\|^p=\left(\sum_{m=0}^p \binom p m \left(\left\| \varphi^{(m)}\right\|_{L^2(\mathbb{R})} \right)^2\right)^{1/2},\quad\varphi\in H^p,
\end{equation*}
and $H^{-p}=\left(H^p\right)^\ast$ with the  norm associated with the strong topology of the adjoint  space.
We have $H^0=L^2(\mathbb{R})=\left(H^0\right)
^\ast=H^{-0}$. By $\langle f,\varphi\rangle$, denote the value of a distribution $f\in H^{-p}$ on a test function $\varphi\in H^p$, $p=0,1$.

By $\widehat{H}^l$, denote the subspace of all even distributions in $H^l$, $l=\overline{-1,1}$. It is easy to see that $\widehat{H}^l$ is a closed subspace of $H^l$, $l=\overline{-1,1}$. 
%Here and further we use the notation $s=\overline{k,m}$ for $s\in \mathbb Z$ such that $k\leq s\leq m$ if $k\leq m$ and for $s\in\emptyset$ if $k>m$, $k\in\mathbb Z$, $m\in\mathbb Z$.

Let $\varphi \in L_{\text{loc}}^2(\Omega)$. We define the  derivative $\mathbb D_{\rho k}$ by the rule
$$
\mathbb D_{\rho k} \varphi = \sqrt{\frac k\rho} \varphi' + Q_1(\rho,k) \varphi.
$$
If, in addition, $\mathbb D_{\rho k} \varphi \in L_{\text{loc}}^2(\Omega)$ and $(\mathbb D_{\rho k} \varphi)'\in   L_{\text{loc}}^2(\Omega)$ (the derivative $(\cdot)'$ is considered in $\EuScript D'(\Omega)$), we can consider  $\mathbb D_{\rho k}^2 \varphi$. Then $\varphi''\in \EuScript D'(\Omega)$ and
$$
\mathbb D_{\rho k}^2 \varphi = \frac1\rho\big(k\varphi'\big)' + Q_2(\rho,k)\varphi.
$$
Obviously, $\mathbb D_{\rho k}^m\varphi= \varphi^{(m)}$ if $\rho=k=1$, $m=0,1$.

Denote 
$$
L_\rho^2(\Omega)=\{f\in L_{\text{loc}}^2(\Omega)\mid \sqrt\rho f\in L^2(\Omega)\} 
$$
with the norm
$$
\| f\|_{L_\rho^2(\Omega)}=\|\sqrt\rho f\|_{L^2(\Omega)}
=\left(\int_\Omega |f(x)|^2\rho(x)\, dx\right)^{1/2},\quad f\in L_\rho^2(\Omega).
$$

For $p=0,1$, consider also the space
\begin{align*}
\overset{\circ}{\mathbb H}{}^p =\{\varphi \in L_{\text{loc}}^2(0,+\infty) &\mid (\forall m= \overline{0,p}\ \mathbb D_{\rho k}^m \varphi \in L_\rho^2(0,+\infty))\\
&\text{and}\ (\forall m= \overline{0,p-1}\ (\mathbb D_{\rho k}^m \varphi)(0^+)=0)\}
\end{align*}
with the norm
$$
\fnl \varphi \fnr^{p \circ}=\left(\sum_{m=0}^p \binom p m \left(\| \mathbb D_{\rho k}^m\varphi\|_{L_\rho^2(\Omega)} \right)^2 \right)^{1/2},\quad \varphi \in \overset{\circ}{\mathbb H}{}^p,
$$
and the dual space $\overset{\circ}{\mathbb H}{}^{-p}= \Big(\overset{\circ}{\mathbb H}{}^p\Big)^*$ with the norm associated with the strong topology of the adjoint space. Evidently, $\overset{\circ}{\mathbb H}{}^0=\overset{\circ}{\mathbb H}{}^{-0}=L_\rho^2(0,+\infty)$. By $\al g, \psi\ar^\circ$, denote the value of a distribution $g\in \overset{\circ}{\mathbb H}{}^{-p}$ on a test function $\psi \in \overset{\circ}{\mathbb H}{}^p$, $p=0,1$. In particular, we have 
$$
\al g, \psi\ar^\circ = \langle g,\psi\rangle_{L_\rho^2(0,+\infty)}
=\int_0^\infty g(x) \psi(x) \rho(x)\, dx,\quad g\in \overset{\circ}{\mathbb H}{}^0,\ \psi \in \overset{\circ}{\mathbb H}{}^0.
$$
Put
$$
\al \mathbb D_{\rho k}f, \varphi\ar^\circ = - \al f, \mathbb D_{\rho k}\varphi\ar^\circ,\quad f\in \overset{\circ}{\mathbb H}{}^0,\ \varphi \in \overset{\circ}{\mathbb H}{}^1.
$$

Consider also the following modified Sobolev spaces introduced and studied in \cite{F1,disF,F3}. Denote
\begin{equation*}
\mathbb{H}^p=\left\lbrace \varphi\in L_\rho^2(\mathbb{R})\mid\forall m=\overline{0,p}\ \mathbb D_{\rho k}^m\varphi\in L_\rho^2(\mathbb{R})\right\rbrace,\quad p=0,1,
\end{equation*}
with the norm
\begin{equation*}
\fnl\varphi\fnr^p=\left(
\sum_{m=0}^p\binom p m\left(\left\|\mathbb D_{\rho k}^m\varphi\right\|_{L_\rho^2(\mathbb{R})}\right)^2
\right)^{1/2},\quad\varphi\in\mathbb{H}^p,\quad p=0,1, 
\end{equation*}
and the dual space $\mathbb{H}^{-p}=\left(\mathbb H^p\right)^\ast$ with the norm associated with the strong topology of the adjoint space.
By $\al f,\varphi\ar$, denote the value of a distribution $f\in\mathbb{H}^{-p}$ on a test function $\varphi\in\mathbb{H}^p$, $p=0,1$. Evidently, $\mathbb{H}^0=\left(\mathbb{H}^0\right)^\ast=L^2(\mathbb{R})$ and
\begin{equation*}
\al f,\varphi\ar=\left\langle f,\varphi\right\rangle_{L_\rho^2(\mathbb{R})}=\int_{-\infty}^\infty f(x)\varphi(x)\rho(x)\,dx,\quad f\in\mathbb{H}^0,\ \varphi\in\mathbb{H}^0.
\end{equation*}
Put 
$$
\al \mathbb D_{\rho k}f,\varphi\ar=-\al f,\mathbb D_{\rho k}\varphi\ar,\quad f\in\mathbb{H}^0,\ \varphi\in\mathbb{H}^1. 
$$
Note that $\mathbb H^{-0}=\mathbb H^0=L_\rho^2(\mathbb{R})$. For $\rho=k=1$, we have $\mathbb{H}^m=H^m$, $m=\overline{-1,1}$.
In \cite{F1}, it has been proved that 
$\mathbb{H}^m\subset\mathbb{H}^n$ is dense continuous embedding, $-1\leq n\leq m\leq 1$, and $\EuScript{D}\subset\mathbb{H}^p\subset\mathbb{H}^{-p}\subset \EuScript{D}'$ are dense continuous embeddings, $p=0,1$, where $\EuScript{D}=\EuScript{D}(\mathbb R)$. However, the relation between the Schwartz space $\EuScript S$ and $\mathbb H^p$ essentially depends on $\rho$ and $k$. For example, if $\rho=k$ then 
$$
\varphi \in \mathbb H^p \Leftrightarrow  \sqrt\rho \varphi\in H^p, \quad p=\overline{-1,1}.
$$
If $\rho(x)=k(x)=\cosh x$, $x\in \mathbb R$, then
$$
\EuScript S\not\subset \mathbb H^p\quad \text{and}\quad \mathbb H^{-p}\not\subset \EuScript S',\quad p=0,1.
$$
If $\rho(x)=k(x)=1/\cosh x$, $x\in \mathbb R$, then 
$$
\EuScript S \subset \mathbb H^p\quad \text{and}\quad \mathbb H^{-p} \subset \EuScript S',\quad p=0,1.
$$

By $\widehat{\mathbb{H}}^s$, denote the subspace of all even distributions in $\mathbb{H}^s$, $s=\overline{-1,1}$. The even extension of a function from $\overset{\circ}{\mathbb H}{}^s$ belongs to $\widehat{\mathbb{H}}^s$, $s=0,1$ (see \cite{F3}). The restriction of a function from $\widehat{\mathbb{H}}^0$ to $[0,+\infty)$ belongs to $\overset{\circ}{\mathbb H}{}^0$. However, there exist functions from $\widehat{\mathbb{H}}^1$ whose restrictions do not belong to  $\overset{\circ}{\mathbb H}{}^1$. Therefore, there exist distributions from $\overset{\circ}{\mathbb H}{}^{-1}$ which cannot be extended to the space $\widehat{\mathbb{H}}^{1}$. But due to the following important theorem proved in \cite[Theorem 3.12]{disF}, the distribution generated by the derivative $\mathbb D_{\rho k}^2 f_+\in\overset{\circ}{\mathbb H}{}^{-1}$  of a function $f_+\in\overset{\circ}{\mathbb H}{}^1$ can be extended to the space $\widehat{\mathbb{H}}^{1}$.
\begin{theorem}
\label{th1ne}
Let $f_+\in\overset{\circ}{\mathbb H}{}^1$, $\varphi\in\widehat{\mathbb H}^1$  and $f$ be the even extension of $f_+$.
If  $\left(D_{\rho k}f_+\right)(0^+)\in\mathbb{R}$, then the distribution $\mathbb D_{\rho k}^2 f_+\in\overset{\circ}{\mathbb H}{}^{-1}$ can be extended to the even distribution $F\in\widehat{\mathbb H}^{-1}$ such that
\begin{equation*}
\al F,\varphi \ar=\al\mathbb D_{\rho k}^2 f,\varphi\ar+2\sqrt{(\rho k)(0)}\left(\mathbb D_{\rho k}f_+\right)(0^+)\varphi(0).
\end{equation*}
\end{theorem}

Put
\begin{equation}
\label{qqqq}
q=Q_2(\rho,  k)-\gamma. 
\end{equation}
Due to \eqref{pot1}, $q\in L^{\infty}(0,+\infty)\cap C^1[0,+\infty)$. Note that $q$ is defined on $\mathbb R$ and $q\in C^1(-\infty,0]\cup C^1[0,+\infty)$, but $q'$ may have a jump at $x=0$.

%Let $w^T\in \overset{\circ}{\mathbb H}{}^1$. Consider the steering condition
%\begin{equation}
%w(\cdot,T)=w^T,\quad x\in(0,+\infty).\label{ec}

We will use the transformation operator $\widehat{\mathbb T}
=\mathbf S \widehat{\mathbf T}_r:\widehat{H}^{-1}
\to\widehat{\mathbb{H}}^{-1}$ 
to investigate controllability problems for system \eqref{eq}--\eqref{ic}. The operators $\mathbf S$ and $\widehat{\mathbf T}_r$ have been introduced and studied in \cite{disF,F3}.
\begin{theorem}[\kern-0.8ex \cite{disF,F3}]
The following assertions hold.
\label{th-trans-neu}
\begin{enumerate}[label=\upshape{(\roman*)},ref=\thetheorem\ \upshape{(\roman*)},leftmargin=6ex, labelwidth=6ex, labelsep=0ex, align=left]
\item 
\label{tt1ne} The operator $\widehat{\mathbb{T}}$ is an isomorphism of $\widehat{H}^m$ and $\widehat{\mathbb{H}}^m$, $m=-1,0,1$.
\item
\label{tt2ne}  $\widehat{\mathbb{T}}\delta=\sqrt[4]{(\rho k)(0)}\delta$.
\item
\label{tt3ne}  If $g\in\widehat{H}^1$ and $g'(0^+)\in \mathbb{R}$, then $\left(\mathcal{D}_{\rho k}\widehat{\mathbb{T}}g\right)(0^+)\in \mathbb{R}$ and
\begin{equation*}
\left(\mathcal{D}^2_{\rho k}-q\right)\widehat{\mathbb{T}}g
-2\sqrt{(\rho k)(0)}\left(\mathcal{D}_{\rho k}\widehat{\mathbb{T}}g\right)(0^+)\delta
=\widehat{\mathbb{T}}\left(\frac{d^2}{d\xi^2}g-2g'(0^+)\delta\right).
\end{equation*}
\item
\label{tt4ne}  If $f\in\widehat{\mathbb{H}}^1$ and $\left(\mathcal{D}_{\rho k}f\right)(0^+)\in \mathbb{R}$, then $\left(\widehat{\mathbb{T}}^{-1} f\right)'(0^+)\in \mathbb{R}$ and
\begin{align*} 
\frac{d^2}{d\xi^2}\widehat{\mathbb{T}}^{-1}f-2\left(\widehat{\mathbb{T}}^{-1}f\right)'&(0^+)\delta
\\
&=\widehat{\mathbb{T}}^{-1}\left(\left(\mathcal{D}^2_{\rho k}-q\right)f-2\sqrt{(\rho k)(0)}\left(\mathcal{D}_{\rho k}f\right)(0^+)\delta\right).
\end{align*}
\end{enumerate}
Here $\delta$ is the Dirac distribution.
\end{theorem}
A description and some properties of the operators $\mathbf S$ and $\widehat{\mathbf T}_r$ are given in Section \ref{operators}.

\subsection{Main results.}\label{neum}
Consider control system \eqref{eq}--\eqref{ic}. We suppose that $\left(\frac d{dt}\right)^p w: [0,T]\to \overset{\circ}{\mathbb H}{}^{1-2p}$, $p=0,1$; $w^0\in \overset{\circ}{\mathbb H}{}^1$. One can easily see that equation \eqref{eq} can be rewritten in the form
\begin{equation}
\label{eqq}
w_t=\mathbb D_{\rho k}^2 w -q w,\quad t\in(0,T),
\end{equation}
and condition \eqref{bc} is equivalent to the condition
\begin{equation}
\label{bcne} 
\left(\mathcal{D}_{\rho k}w\right)(0,\cdot)=u,\quad t\in(0,T).
\end{equation}
Let $w^T\in \overset{\circ}{\mathbb H}{}^1$. Consider the steering condition for system \eqref{eq}--\eqref{ic}
\begin{equation}
\label{ec}
w(\cdot,T)=w^T, \quad x\in(0,+\infty).
\end{equation}

Let $w(\cdot,t),w^0\in\overset{\circ}{\mathbb H}{}^1$ and let $W(\cdot,t), W^0$ be their even extensions with respect to $x$, respectively, $t\in [0,T]$. Let $q$ be defined by \eqref{qqqq}. If $w$ is a solution to control system \eqref{eq}--\eqref{ic}, then using Theorem \ref{th1ne} and taking into account \eqref{eqq} and  \eqref{bcne}, we conclude that $W$ is a solution to the system
\begin{align}
&W_t=\mathcal{D}^2_{\rho k}W-qW-2\sqrt{(\rho k)(0)}u\delta,&& \text{on }\mathbb{R}\times(0,T),\label{eq1ne}\\
&W(\cdot,0)=W^0, && \text{on }\mathbb{R},\label{ic1ne}
\end{align}
where $\left(\frac d{dt}\right)^l W:[0,T]\to\widehat{\mathbb{H}}^{1-2l}$, $l=0,1$, 
$W^0\in\widehat{\mathbb{H}}^1$, $\delta$ is the Dirac distribution with respect to $x$.
Let $W (\cdot,t),W^0\in\widehat{\mathbb{H}}^1$ and let $w(\cdot,t),w^0$ be their restrictions to $(0,+\infty)$ with respect to $x$, respectively, $t\in [0,T]$. If $W$ is a solution to control system \eqref{eq1ne}, \eqref{ic1ne}, then due to Corollary \ref{corol1ne} (see Section \ref{operators} below),
\begin{equation}
\label{bc1ne}
\left(\mathcal{D}_{\rho k}w\right)(0,\cdot)=\left(\mathcal{D}_{\rho k}W\right)(0^+,\cdot)=u \quad\text{a.e. on } (0,T).
\end{equation}
Hence, $w$ is a solution to control system \eqref{eq}--\eqref{ic}. 

Let $w^T\in \overset{\circ}{\mathbb H}{}^1$ and let $W^T\in\widehat{\mathbb{H}}^1$ be its  even extension with respect to $x$. It is easy to see that $w(\cdot,T)=w^T$ iff $W(\cdot,T)=W^T$.

Thus, control systems \eqref{eq}--\eqref{ic} and \eqref{eq1ne}, \eqref{ic1ne} are equivalent. Taking into account this equivalence, we will further consider  system \eqref{eq1ne}, \eqref{ic1ne}.

Let $T>0$, $W^0\in\widehat{\mathbb{H}}^1$. By $\EuScript{R}_T^{\rho k\gamma}(W^0)$, denote the set of all states $W^T\in\widehat{\mathbb{H}}^1$ for which there exists a control $u^{\rho k\gamma}\in L^\infty(0,T)$ such that there exists a unique solution $W$ to \eqref{eq1ne}, \eqref{ic1ne} with $u=u^{\rho k\gamma}$ and $W(\cdot,T)=W^T$.

\begin{definition}
A state $W^0\in\widehat{\mathbb{H}}^1$ is said to be controllable to a state $W^T\in\widehat{\mathbb{H}}^1$ with respect to system \eqref{eq1ne}, \eqref{ic1ne} in a given time $T>0$ if $W^T\in\EuScript{R}_T^{\rho k \gamma}(W^0)$.
\end{definition}

\begin{definition}
A state $W^0\in\widehat{\mathbb{H}}^1$ is said to be approximately controllable to a state $W^T\in\widehat{\mathbb{H}}^1$ with respect to system \eqref{eq1ne}, \eqref{ic1ne} in a given time 
$T>0$ if $W^T\in\overline{\EuScript{R}_T^{\rho k\gamma}(W^0)}$, where the closure is considered in the space $\widehat{\mathbb{H}}^1$.
\end{definition}
Thus, the main goal of the paper is to investigate whether the state $W^0$ is  controllable (approximately controllable) to a target state $W^T$ with respect to system \eqref{eq1ne}, \eqref{ic1ne} in a given time $T$.

To this aid, consider the  control system with the simplest heat operator (system \eqref{eq1ne}, \eqref{ic1ne} with $\rho=k=1$, $\gamma=0$): 
\begin{align}
&Z_t=Z_{\xi\xi}-2u\delta,&&\text{on }\mathbb{R}\times(0,T),\label{z1ne}\\
&Z(\cdot,0)=Z^0, &&\text{on }\mathbb{R},\label{z2ne}
\end{align}
where $u\in L^\infty(0,T)$ is a control, $u=u^{110}$, $\left(\frac{d}{dt}\right)^l Z:[0,T]\rightarrow \widehat{H}^{1-2l}$, $l=0,1$, $Z^0\in\widehat{H}^1$. Let $Z^T\in\widehat{H}^1$. Consider also the steering condition for this system:
\begin{equation*}
Z(\cdot,T)=Z^T, \quad\text{on }\mathbb{R}.
\end{equation*}
Control  system \eqref{z1ne}, \eqref{z2ne} has been  investigated in \cite{FKh2}.  In particular, it has been proved therein that
\begin{equation}
 Z_x(0^+,\cdot)=u,\quad \text{a.e. on } (0,T).\label{z4ne}\\
\end{equation}

Using Theorems \ref{thtransf1ne} and \ref{thtransf2ne}  (see Section \ref{operators} below), we obtain the following theorem.
\begin{theorem}
\label{corol2ne}
Let $T>0$, $W^0\in\widehat{\mathbb{H}}^1$, $W^T\in\widehat{\mathbb{H}}^1$, $Z^0=\widehat{\mathbb{T}}^{-1}W^0$, $Z^T=\widehat{\mathbb{T}}^{-1}W^T$. Then
\begin{enumerate}[label=\upshape{(\roman*)},ref=\upshape{(\roman*)},leftmargin=6ex, labelwidth=6ex, labelsep=0ex, align=left]
\item\label{daa0ne}
$\EuScript{R}_T^{\rho k\gamma}\left(W^0\right)=\widehat{\mathbb{T}}\left(\EuScript{R}_T^{110}\left(Z^0\right)\right)$.
\item\label{daa1ne}
A state $Z^0$ is controllable to a state $Z^T$ with respect to system \eqref{z1ne}, \eqref{z2ne} in a time $T$ iff a state $W^0$ is controllable to a state $W^T$ with respect to system \eqref{eq1ne}, \eqref{ic1ne} in this time $T$.
\item\label{daa3ne}
A state $Z^0$ is approximately controllable to a state $Z^T$ with respect to system \eqref{z1ne}, \eqref{z2ne} in a time $T$ iff a state $W^0$ is approximately controllable to a state $W^T$ with respect to system \eqref{eq1ne}, \eqref{ic1ne} in this time $T$.
\end{enumerate}
\end{theorem}

Thus, control system \eqref{eq1ne}, \eqref{ic1ne} with a general heat operator replicates the controllability properties of control system \eqref{z1ne}, \eqref{z2ne} with the simplest heat operator and vice versa. 

The main results of the paper are the following two theorems.
\begin{theorem}
\label{lackne}
If a state $W^0\in\widehat{\mathbb{H}}^1$ is controllable to $0$ with respect to system \eqref{eq1ne}, \eqref{ic1ne} in a time $T>0$, then $W^0=0$.
\end{theorem}
\begin{theorem}
\label{thapprne}
Each state $W^0\in\widehat{\mathbb{H}}^1$ is  approximately controllable  to any target state $W^T\in\widehat{\mathbb{H}}^1$ with respect to system \eqref{eq1ne}, \eqref{ic1ne} in a given time $T>0$.
\end{theorem}

In the case $\rho=k=1$, $\gamma=0$ these theorems have been proved in \cite{FKh2}.
By using Theorem  \ref{corol2ne}, we obtain Theorems \ref{lackne} and \ref{thapprne}.

Taking into account the algorithm given in \cite[Section 7]{FKh2}, one can construct  piecewise constant controls solving   the  approximate controllability problem for system \eqref{z1ne}, \eqref{z2ne}. Hence, using Theorem \ref{thtransf1ne}, one can obtain  controls solving the  approximate controllability problem for system  \eqref{eq1ne}, \eqref{ic1ne} (see Section \ref{operators} below).

\section{The transformation operator \texorpdfstring{$\widehat{\mathbb{T}}$}{hat T} and it's application to a control system}\label{operators}

In this section, we recall some properties of the operator $\widehat{\mathbb T}$ and apply it to control system \eqref{eq1ne}, \eqref{ic1ne}.
We have $\widehat{\mathbb T}
=\mathbf S \widehat{\mathbf T}_r:\widehat{H}^{-1}
\to\widehat{\mathbb{H}}^{-1}$.

The operator $\mathbf S: H^{-1}\to \mathbb H^{-1}$ has been introduced and studied in \cite{disF,F3}.

\begin{theorem}[\kern-0.8ex\cite{disF,F3}] 
\label{sss}
The following assertions hold.
\begin{enumerate}[label=\upshape{(\roman*)}, ref=\thetheorem\ \upshape{(\roman*)},leftmargin=6ex, labelwidth=6ex, labelsep=0ex, align=left]
\item \label{sss-i}
The operator $\mathbf S$ is an isometric isomorphism of $H^m$ and $\mathbb H^m$, $m=\overline{-1,1}$;
\item \label{sss-ii}
$ \mathcal{D}_{\rho k}\mathbf S \psi=\mathbf S \frac d{d\lambda}\psi$, $\psi \in H^m$, $m=0,1$;
\item \label{sss-iii}
$\al f, \varphi\ar=\left\langle \mathbf S^{-1}f, \mathbf S^{-1}\varphi \right\rangle$, $f \in \mathbb H^{-m}$, $\varphi \in \mathbb H^m$, $m=0,1$;
\item \label{sss-iv}
$\mathbf S \delta= \sqrt[4]{(\rho k)(0)}\delta$.
\end{enumerate}
\end{theorem}

In particular, we have
\begin{equation*}
\mathbf S \psi=\frac{\psi\circ\sigma}{\sqrt[4]{\rho k}},\quad \psi\in H^0,\quad
\text{and}\quad
\mathbf S^{-1} \varphi=\left(\sqrt[4]{\rho k}\varphi\right)\circ\sigma^{-1},\quad \varphi\in\mathbb{H}^0,
\end{equation*}
where $\psi\circ\sigma=\psi(\sigma(x))$, $\sigma$ is defined by \eqref{sigma}. It follows from  \eqref{sigma}, \eqref{sigma-a} that $\sigma$ is an odd increasing invertible function and $\sigma(x)\rightarrow\pm\infty$ as $x\rightarrow\pm\infty$.

Put
\begin{equation}
\label{r}
r(\lambda)=\left(q\circ\sigma^{-1}\right)(\lambda)=\left(\left(Q_2( \rho,k)-\gamma\right)\circ\sigma^{-1}\right)(\lambda),\quad\lambda\in[0,+\infty).
\end{equation}
Due to \eqref{pot1} and \eqref{pot2}, we have
\begin{equation}
\label{rr}
r\in L^{\infty}(0,+\infty)\cap C^1[0,+\infty)\quad\text{and}\quad\lambda r\in L^1(0,+\infty).
\end{equation}
Consider the operator $\widehat{\mathbf T}_r: \widehat H^{-1}\to \widehat H^{-1}$. This operator is the extension to $\widehat{H}^{-1}$ of the well-known transformation operator of the Sturm--Liouville problem (see, e.g.,  \cite[Chap. 3]{M}).
The complete description of the extension and its application to the wave equation with variable coefficients have been given in \cite{F3,disF,Kh2}.

\begin{theorem}[\kern-0.8ex\cite{F3,disF}] 
\label{trans-neu}
The following assertions hold.
\begin{enumerate}[label=\upshape{(\roman*)},ref=\thetheorem\ \upshape{(\roman*)},leftmargin=6ex, labelwidth=6ex, labelsep=0ex, align=left]
\item \label{trans-neu-i}
The operator $\widehat{\mathbf T}_r$ is an automorphism of $\widehat H^m$, $m=\overline{-1,1}$.
\item \label{trans-neu-ii}
 If $g\in\widehat{H}^1$ and $g'(0^+)\in \mathbb{R}$, then $\left(\widehat{\mathbf T}_rg\right)'(0^+)\in \mathbb{R}$ and
\begin{equation*}
\left(\frac{ d^2}{d\lambda^2}-r\right)\widehat{\mathbf{T}}_rg
-2\left(\widehat{\mathbf T}_rg\right)'(0^+)\delta
=\widehat{\mathbf T}_r\left(\frac{d^2}{d\xi^2}g-2g'(0^+)\delta\right).
\end{equation*}
\item \label{trans-neu-iii}
If $f\in\widehat{H}^1$ and $f'(0^+)\in \mathbb{R}$, then $\left(\widehat{\mathbf T}_r^{-1}f\right)'(0^+)\in \mathbb{R}$ and
\begin{equation*}
\frac{d^2}{d\xi^2}\widehat{\mathbf T}_r^{-1}f-2\left(\widehat{\mathbf T}_r^{-1}f\right)'(0^+)\delta=\widehat{\mathbf T}_r^{-1}\left(\left(\frac{ d^2}{d\lambda^2}-r\right)f-2f'(0^+)\delta\right).
\end{equation*}
\item \label{trans-neu-iv}
$\widehat{\mathbf T}_r\delta=\delta$.
\end{enumerate}
\end{theorem}

In particular, we have
\begin{align*}
\left(\widehat{\textbf{T}}_r g\right)(\lambda)&=g(\lambda)+ \int_{|\lambda|}^\infty K(|\lambda|,\xi)
g(\xi)d\xi,&&\lambda\in\mathbb{R},\ g\in \widehat{H}^0,\\
\left(\widehat{\textbf{T}}_r^{-1} f\right)(\xi)&=f(\xi)+ \int_{|\xi|}^\infty L(|\xi|,\lambda)
f(\lambda)d\lambda,&&\xi\in\mathbb{R},\ f\in \widehat{H}^0,
\end{align*}
where, according to \cite[Chap. 3]{M},  the kernel $K\in C^2(\Omega)$ is a unique solution to the system
\begin{equation}
\label{kern}
\left\lbrace 
\begin{aligned}
&K_{y_1y_1}-K_{y_2y_2}=r(y_1)K,  &&\text{on }\Omega,\\ 
&K(y_1,y_1)=\frac{1}{2}\int_{y_1}^\infty r(\xi)d\xi, && y_1>0,\\
&\lim\limits_{y_1+y_2\rightarrow\infty}K_{y_1}(y_1,y_2)=\lim\limits_{y_1+y_2\rightarrow\infty}K_{y_2}(y_1,y_2)=0, &&\text{on }\Omega,
\end{aligned}
\right.
\end{equation}
 $\Omega=\{(y_1,y_2)\in\mathbb{R}^2\mid y_2>y_1>0\}$, and the kernel $L\in C^2(\Omega)$ is determined by the following equation
\begin{align}
L(y_1,y_2)+K(y_1,y_2)+\int_{y_1}^{y_2} L(y_1,\xi)K(\xi,y_2)d\xi&=0, &&\text{on }\Omega.\label{lk1}
\end{align}
We also need  the following estimates proved in \cite[Chap.~3]{M}:
\begin{align}
|K(y_1,y_2)|&\leq M_0\sigma_0\left(\frac{y_1+y_2}{2} \right),&&\text{on }\Omega,\label{estk1}\\
|K_{y_1}(y_1,y_2)&|\leq\frac{1}{4}\left|r\left(\frac{y_1+y_2}{2} \right)\right|+ M_1\sigma_0\left(\frac{y_1+y_2}{2} \right),&&\text{on }\Omega,\label{estk1ne}
\end{align}
where $M_0>0$, $M_1>0$ are constants, and 
\begin{equation}
\label{sigmnul}
\sigma_0(x)=\int_x^\infty |r(\xi)| d\xi, \quad x>0.
\end{equation}

In the following theorems, the application of the transformation operator $\widehat{\mathbb{T}}$ to a control system is considered.

%%%%%%%%%%%% THEOREM
\begin{theorem}
\label{thtransf1ne}
Let $Z$ be a solution to \eqref{z1ne}, \eqref{z2ne} with $u=u^{110}$, where $u^{110}\in L^\infty(0,T)$, $Z^0\in\widehat{H}^1$. Let $W(\cdot,t)=\left(\widehat{\mathbb{T}}Z\right)(\cdot,t)$, $t\in[0,T]$, $W^0=\widehat{\mathbb{T}}Z^0$. Then $W$ is a solution to system \eqref{eq1ne}, \eqref{ic1ne} with   the control $u=u^{\rho k\gamma}$,
\begin{align}
u^{\rho k\gamma}(t)=\frac{1}{\sqrt[4]{(\rho k)(0)}}&\left(u^{110}(t)+\int_0^\infty K_{y_1}(0,x)Z(x,t)dx\right.\nonumber\\
&-\left.\frac{1}{2}Z(0^+,t)\int_0^\infty r(\xi)d\xi\right),\quad t\in[0,T],\label{contr1ne}
\end{align}
where $K$ is a solution to \eqref{kern}, $r$ is defined by \eqref{r}.
Besides, \eqref{bc1ne} holds and
\begin{align}
\fnl W(\cdot,t)\fnr^1&\leq E_0\|Z(\cdot,t)\|^1,\quad t\in[0,T],\label{est1ne}\\
\|u^{\rho k\gamma}\|_{L^\infty(0,T)}&\leq G_0(T)\|u^{110}\|_{L^\infty(0,T)}+E_1\left\|Z^0\right\|^1,\label{est2ne}
\end{align}
where $E_0 >0$ and $E_1 >0$ are constants independent of $T$, 
$$
G_0(T)=\frac{1}{\sqrt[4]{(\rho k)(0)}}\left(1+(T+3)\left(\frac{2\sqrt{\sigma_0(0)}}{\sqrt{\pi}}\sqrt{R_0+M_1^2R}+\frac{\sigma_0(0)}{\sqrt{2\pi}}\right)\right),
$$
$M_1$ is the constant from \eqref{estk1ne}, $\sigma_0$ is defined by \eqref{sigmnul}, and
\begin{equation}
\label{rconst}
R=\int_0^\infty\xi |r(\xi)|d\xi,\quad
R_0=\frac{1}{16}\|r\|_{L^\infty(0,+\infty)}.
%\label{rconstin}
\end{equation}
\end{theorem}

\begin{proof}%[Proof of Theorem \ref{thtransf1ne}]
The first part of this theorem is proved similarly to the first part of the corresponding theorem in \cite{F3,disF} (\cite[Theorem 6.12]{disF}, \cite[Theorem 4.2]{F3}).   
Applying Theorem \ref{tt3ne}, we obtain the first assertion of the theorem.

Taking into account Theorem \ref{sss-ii}, \eqref{z4ne}, \eqref{kern}, and \eqref{contr1ne} we obtain
\begin{align*}
\left(\mathcal{D}_{\rho k}W\right)(0^+,t)=\left(\mathcal{D}_{\rho k}\widehat{\mathbb{T}}Z\right)(0^+,t)=\textbf{S}\left(\widehat{\textbf{T}}_rZ\right)'(0^+,t)=\frac{1}{\sqrt[4]{(\rho k)(0)}}\left(\vphantom{\int_0^\infty}Z_x(0^+,t)\right.\\
+\left.\int_0^\infty K_{y_1}(0,x)Z(x,t)dx-K(0,0)Z(0^+,t)\right)=u^{\rho k\gamma}(t),\quad t\in[0,T].
\end{align*}
Thus, \eqref{bc1ne} is valid.

It follows from Theorem \ref{tt1ne} that there exists a constant $E_0>0$ such that \eqref{est1ne} holds.

To complete the proof, it remains to prove  \eqref{est2ne}. Due to \eqref{estk1ne},  we obtain from \eqref{contr1ne}
\begin{align*}
&\|u^{\rho k\gamma}\|_{L^\infty(0,T)}\leq\frac{1}{\sqrt[4]{(\rho k)(0)}}\left(\vphantom{\sqrt{\int_0^\infty\left|\frac{1}{4}r\left(\frac{x}{2}\right)+ M_1\sigma_0\left(\frac{x}{2}\right)\right|^2 dx}}\|u^{110}\|_{L^\infty(0,T)}\right.\\
&+\left.\frac{\|Z(\cdot,t)\|^0}{\sqrt{2}}\sqrt{\int_0^\infty\left|\frac{1}{4}r\left(\frac{x}{2}\right)+ M_1\sigma_0\left(\frac{x}{2}\right)\right|^2 dx}
+\frac{1}{2}\sigma_0(0)\left|Z(0^+,t)\right|\right),\ t\in[0,T].     
\end{align*}
Since $\|Z(\cdot,t)\|^0\leq \|Z(\cdot,t)\|^1$ (\cite[Chap. 1]{VG}), $t\in[0,T]$, we get from here that
\begin{align}
\|u^{\rho k\gamma}\|_{L^\infty(0,T)}\leq\frac{1}{\sqrt[4]{(\rho k)(0)}}&\left(\|u^{110}\|_{L^\infty(0,T)}+\|Z(\cdot.t)\|^1\sqrt{R_0\sigma_0(0)+M_1^2\sigma_0(0)R}\right.\nonumber\\
&+\left.\frac{1}{2}\sigma_0(0)|Z(0^+,t)|\right). \label{esne11}    
\end{align}
For  $Z\in\widehat{H}^1$ we have
$Z(0^+,t)=\frac{1}{\sqrt{2\pi}}\int_{-\infty}^\infty(\EuScript{F}Z)(\sigma,t)d\sigma$, $t\in[0,T]$, where $\EuScript F:H^0\to H^0$ is the Fourier transform operator, and $\EuScript F H^1=H_1$, $H_1=\{f\in H^0\mid (1+|\sigma|^2)^{1/2}f\in H^0\}$, $\|f\|_1=\|(1+|\sigma|^2)^{1/2}f\|^0$, $f\in H_1$ (see, e.g., \cite[Chap.~1]{VG}).
Hence,
\begin{align}
&|Z(0^+,t)|=\frac{1}{\sqrt{2\pi}}\left|\int_{-\infty}^\infty\sqrt{1+\sigma^2}(\EuScript{F}Z)(\sigma,t)\frac{d\sigma}{\sqrt{1+\sigma^2}}\right|\nonumber\\
&\leq\frac{1}{\sqrt{2\pi}}\|(\EuScript{F}Z)(\cdot,t)\|_1\sqrt{\int_{-\infty}^\infty\frac{d\sigma}{1+\sigma^2}}=\frac{1}{\sqrt{2}}\|Z(\cdot,t)\|^1,\quad t\in[0,T].\label{esmod}
\end{align}
Substituting \eqref{esmod} into \eqref{esne11}, we get
\begin{align}
&\|u^{\rho k\gamma}\|_{L^\infty(0,T)}
\leq\frac{1}{\sqrt[4]{(\rho k)(0)}}\left(\vphantom{\frac{\sigma_0(0)}{2\sqrt{2}}}\|u^{110}\|_{L^\infty(0,T)}\right.\nonumber\\
&+\left.\|Z(\cdot,t)\|^1\left(\sqrt{\sigma_0(0)}\sqrt{R_0+M_1^2R}
+\frac{\sigma_0(0)}{2\sqrt{2}}\right)\right),\quad t\in[0,T].\label{esne}
\end{align}
Using formula (13) from \cite{FKh2}, we have
\begin{equation*}
(\EuScript{F}Z)(\sigma,t)=e^{-\sigma^2t}(\EuScript{F}Z^0)(\sigma)-\sqrt{\frac{2}{\pi}}\int_0^t e^{-(t-\xi)\sigma^2}u^{110}(\xi)d\xi,\quad \sigma\in\mathbb{R},\ t\in[0,T].
\end{equation*}
It is easy to obtain from here 
\begin{align}
\|Z(\cdot,t)\|^1&=\|(\EuScript{F}Z)(\cdot,t)\|_1\leq\left\|\EuScript{F}Z^0\right\|_1+\frac{2(T+3)}{\sqrt{\pi}}\|u^{110}\|_{L^\infty(0,T)}\nonumber\\
&=\left\|Z^0\right\|^1+\frac{2(T+3)}{\sqrt{\pi}}\|u^{110}\|_{L^\infty(0,T)},\quad t\in[0,T].\label{estzzne}
\end{align}
Substituting  \eqref{estzzne} into  \eqref{esne}, we get
\begin{align*}
\|u^{\rho k\gamma}\|_{L^\infty(0,T)}\leq\frac{1}{\sqrt[4]{(\rho k)(0)}}&\left(\|u^{110}\|_{L^\infty(0,T)}+\left(\left\|Z^0\right\|^1+\frac{2(T+3)}{\sqrt{\pi}}\|u^{110}\|_{L^\infty(0,T)}\right)\right.\\
&\times\left.\left(\sqrt{\sigma_0(0)(R_0+M_1^2R)}
+\frac{\sigma_0(0)}{2\sqrt{2}}\right)\right).  
\end{align*}
The theorem is proved.
\end{proof}

%%%%%%%%%%%% COROLLARY
\begin{corollary}
\label{corol1ne}
Let $W$ be a solution to \eqref{eq1ne}, \eqref{ic1ne} with $u=u^{\rho k\gamma}$, where  $W^0\in\widehat{\mathbb{H}}^1$ and $u\in L^\infty(0,T)$. Then $\left(\mathcal{D}_{\rho k}W\right)(0^+,\cdot)=u$  a.e. on $[0,T]$, i.e., \eqref{bc1ne} holds.
\end{corollary}

\begin{proof}%[Proof of Corollary \ref{corol1ne}]
Put $Z(\cdot,t)=\left(\widehat{\mathbb{T}}^{-1}W\right)(\cdot,t)$, $t\in[0,T]$ and apply the operator $\widehat{\mathbb{T}}^{-1}$ to \eqref{eq1ne}. Due to Theorem \ref{tt4ne}, we obtain
\begin{align*}
Z_t(\cdot,t)&=Z_{\xi\xi}(\cdot,t)-2Z_x(0^+,t)\delta
\\
&\kern2.ex+2\sqrt{(\rho k)(0)}\left(\left(\mathcal{D}_{\rho k}W\right)(0^+,t)-u^{\rho k\gamma}(t)\right)\widehat{\mathbb{T}}^{-1}\delta,\quad t\in[0,T]. 
\end{align*}
Using Theorem \ref{tt2ne}, we get
\begin{align*} 
Z_t(\cdot,t)&=Z_{\xi\xi}(\cdot,t)-2\left(Z_x(0^+,t)
\right.
\\
&\kern2.8ex\left.-\sqrt[4]{(\rho k)(0)}\left(\mathcal{D}_{\rho k}W\right)(0^+,t)+\sqrt[4]{(\rho k)(0)}u^{\rho k\gamma}(t)\right)\delta,\quad t\in[0,T].
\end{align*}
Thus, $Z$ is a solution to system \eqref{z1ne}, \eqref{z2ne} with $Z^0=\widehat{\mathbb{T}}^{-1}W^0$ and with the control $u=u^{110}$,
\begin{equation*}
u^{110}(t)=Z_x(0^+,t)-\sqrt[4]{(\rho k)(0)}\left(\mathcal{D}_{\rho k}W\right)(0^+,t)+\sqrt[4]{(\rho k)(0)}u^{\rho k\gamma}(t),\quad t\in[0,T].
\end{equation*}
Due to \eqref{z4ne}, we get $\left(\mathcal{D}_{\rho k}W\right)(0^+,t)=u^{\rho k\gamma}(t)$, $t\in[0,T]$.
\end{proof}

To prove the next theorem, we need the following lemma proved in \cite{FKh3}.
\begin{lemma}
\label{exist}
Let 
\begin{equation*}
|f(t)|\leq N_0\quad\text{and}\quad |P(t)|\leq\frac{N_1}{\sqrt{\pi t}},\quad t\in[0,T],
\end{equation*}
where $N_0>0$ and $N_1>0$ are constants.  Then there exists a unique solution $v\in L^\infty(0,T)$ to equation
\begin{equation}
\label{eqoper}
v(t)=f(t)+\int_0^t v(\xi)P(t-\xi)d\xi,\quad t\in[0,T],
\end{equation}
and
\begin{equation}
\label{estvvv}
\|v\|_{L^\infty(0,T)}\leq N_0\left(1+2N_1\sqrt{\frac{T}{\pi}}e^{N_1^2T}\right).
\end{equation}
\end{lemma}

%%%%%%%%%%%%% THEOREM
\begin{theorem}
\label{thtransf2ne}
Let $W$ be a solution to \eqref{eq1ne}, \eqref{ic1ne} with $u=u^{\rho k\gamma}$, where $u^{\rho k\gamma}\in L^\infty(0,T)$, $W^0\in\widehat{\mathbb{H}}^1$. Let $Z(\cdot,t)=\left(\widehat{\mathbb{T}}^{-1}W\right)(\cdot,t)$, $t\in[0,T]$, $Z^0=\widehat{\mathbb{T}}^{-1}W^0$. Then $Z$ is a solution to system \eqref{z1ne}, \eqref{z2ne} with the control $u=u^{110}$,
\begin{align}
u^{110}(t)=\sqrt[4]{(\rho k)(0)}u^{\rho k\gamma}(t)+\frac{1}{2}\sqrt[4]{(\rho k)(0)}W(0^+,t)\int_0^\infty r(\mu)d\mu\nonumber\\
+\int_0^\infty L_{y_1}(0,x)\left(\textbf{S}^{-1}W\right)(x,t)dx,\quad  t\in[0,T],\label{contr2ne}
\end{align}
where $L$ is determined by \eqref{lk1}, $r$ is defined by \eqref{r}. In addition, 
\begin{align}
\|Z(\cdot,t)\|^1&\leq E_2\fnl W(\cdot,t)\fnr^1,\quad t\in[0,T],\label{est3ne}\\
\|u^{110}\|_{L^\infty(0,T)}&\leq G_1(T)\left(\|u^{\rho k\gamma}\|_{L^\infty(0,T)}+E_3\fnl W^0\fnr^1\right),\label{est4ne}
\end{align}
where $E_2 >0$ and $E_3 >0$ are constants independent of $T$, 
$$
G_1(T)=\sqrt[4]{(\rho k)(0)}e^{(\sigma_0(0)+2M_1R)^2T}\left(1+2\sqrt{\frac{T}{\pi}}\big(\sigma_0(0)+2M_1R\big)\right),
$$
$M_1$ is the constant from \eqref{estk1ne}, $\sigma_0$ is defined by \eqref{sigmnul}, $R$ is defined by \eqref{rconst}.
\end{theorem}

\begin{proof}%[Proof of Theorem \ref{thtransf2ne}]
Applying Theorem \ref{th-trans-neu} 
(see \ref{tt1ne}, \ref{tt4ne}), and Corollary \ref{corol1ne} we obtain \eqref{contr2ne} and \eqref{est3ne}. Let us prove \eqref{est4ne}.
From \eqref{contr2ne}, it follows that
\begin{align*}
u^{110}(t)&=\sqrt[4]{(\rho k)(0)}u^{\rho k\gamma}(t)+\frac{1}{2}\left(\widehat{\textbf{T}}_rZ\right)(0^+,t)\int_0^\infty r(\mu)d\mu\\
+&\int_0^\infty L_{y_1}(0,x)\left(\widehat{\textbf{T}}_rZ\right)(x,t)dx=\sqrt[4]{(\rho k)(0)}u^{\rho k\gamma}(t)+\frac{1}{2}Z(0^+,t)\int_0^\infty r(\mu)d\mu\\
+&\frac{1}{2}\int_0^\infty r(\mu)d\mu\int_0^\infty K(0,x)Z(x,t)dx+\int_0^\infty L_{y_1}(0,x)Z(x,t)dx\\
+&\int_0^\infty Z(x,t)\int_0^x L_{y_1}(0,\xi)K(\xi,x)d\xi dx,\quad t\in[0,T].
\end{align*}
By differentiating \eqref{lk1} with respect to $y_1$, we get
\begin{equation*}
-K_{y_1}(0,x)=L_{y_1}(0,x)+\frac{1}{2}K(0,x)\int_0^\infty r(\mu)d\mu +\int_0^x L_{y_1}(0,\xi)K(\xi,x)d\xi,\quad x>0.
\end{equation*}
Therefore,
\begin{align*}
u^{110}(t)=\sqrt[4]{(\rho k)(0)}u^{\rho k\gamma}(t)&+\frac{1}{2}Z(0^+,t)\int_0^\infty r(\mu)d\mu
\\
&-\int_0^\infty K_{y_1}(0,x)Z(x,t)dx,\quad t\in[0,T].
\end{align*}
 (In fact, it is relation \eqref{contr1ne} from Theorem \ref{thtransf1ne}.)

Using formula (15) (for solution to \eqref{z1ne}, \eqref{z2ne}) from \cite{FKh2}, we have
\begin{equation*}
Z(x,t)=\frac{e^{-\frac{x^2}{4t}}}{\sqrt{4\pi t}}\ast Z^0(x)-\sqrt{\frac{2}{\pi}}
\int_0^t u^{110}(\xi)\frac{e^{-\frac{x^2}{4(t-\xi)}}}{\sqrt{2(t-\xi)}}d\xi,\quad x\in\mathbb{R},\ t\in[0,T].
\end{equation*}
Thus, we obtain
\begin{align}
u^{110}(t)&=\sqrt[4]{(\rho k)(0)}u^{\rho k\gamma}(t)+\frac{1}{2}\int_0^\infty r(\mu)d\mu\int_{-\infty}^\infty \frac{e^{-\frac{x^2}{4t}}}{\sqrt{4\pi t}}Z^0(x)dx\nonumber\\
&-\frac{1}{2\sqrt{\pi}}\int_0^\infty r(\mu)d\mu\int_0^t\frac{u^{110}(\xi)}{\sqrt{t-\xi}}d\xi
\nonumber-\int_0^\infty K_{y_1}(0,x)\left(\frac{e^{-\frac{x^2}{4t}}}{\sqrt{4\pi t}}\ast Z^0(x)\right)dx\nonumber\\
&+\sqrt{\frac{2}{\pi}}
\int_0^\infty K_{y_1}(0,x)\int_0^t u^{110}(\xi)\frac{e^{-\frac{x^2}{4(t-\xi)}}}{\sqrt{2(t-\xi)}}d\xi dx,\quad  t\in[0,T].\label{contrneu}
\end{align}
Denote
\begin{align}
f(t)&=\sqrt[4]{(\rho k)(0)}u^{\rho k\gamma}(t)+\frac{1}{2}\int_0^\infty r(\mu)d\mu\int_{-\infty}^\infty \frac{e^{-\frac{x^2}{4t}}}{\sqrt{4\pi t}}Z^0(x)dx
\nonumber\\
&
\kern14.5ex
-\int_0^\infty K_{y_1}(0,x)\left(\frac{e^{-\frac{x^2}{4t}}}{\sqrt{4\pi t}}\ast Z^0(x)\right)dx,&& t\in[0,T],\label{ffneu}\\
P(t)&=\frac{1}{\sqrt{\pi t}}\left(\int_0^\infty K_{y_1}(0,x)e^{-\frac{x^2}{4t}}dx-\frac{1}{2}\int_0^\infty r(\mu)d\mu\right),&& t\in[0,T].\label{ppneu}
\end{align}
Then \eqref{contrneu} takes the form \eqref{eqoper}.
Let us estimate $f$ and $P$. We have
\begin{align}
&\left|\int_{-\infty}^\infty \frac{e^{-\frac{x^2}{4t}}}{\sqrt{4\pi t}}Z^0(x)dx\right|=\frac{1}{\sqrt{2\pi}}\left|\int_{-\infty}^\infty e^{-t\sigma^2}\left(\EuScript{F}Z^0\right)(\sigma)d\sigma\right|\nonumber\\
&\qquad\qquad\leq\frac{1}{\sqrt{2\pi}}\int_{-\infty}^\infty \sqrt{1+\sigma^2}\left|\EuScript{F}Z^0\right|(\sigma)\frac{d\sigma}{\sqrt{1+\sigma^2}}
\nonumber\\
&\qquad\qquad\leq\frac{1}{\sqrt{2\pi}}\left\|\EuScript{F}Z^0\right\|_1\sqrt{\int_{-\infty}^\infty\frac{d\sigma}{1+\sigma^2}}=\frac{1}{\sqrt{2}}\left\|Z^0\right\|^1,\quad t\in[0,T].\label{ttne1}
\end{align}
According to \eqref{estk1ne}, we get
\begin{align}
\left(\left\|K_{y_1}(0,\cdot)\right\|_{L^2(0.+\infty)}\right)^2&\leq\frac{1}{16}\int_0^\infty r^2\left(\frac{x}{2}\right)dx+M_1^2\int_0^\infty \sigma_0^2\left(\frac{x}{2}\right)dx\nonumber\\
&\leq 2\sigma_0(0)\left(R_0+M_1^2R\right),\label{ttne2}
\end{align}
where $R_0$ is defined by \eqref{rconst}. We also have
\begin{equation}
\label{ttne3}
\Bigg\|\frac{e^{-\frac{(\cdot)^2}{4t}}}{\sqrt{4\pi t}}*Z^0\Bigg\|^0
=\left\|e^{-t(\cdot)^2}\EuScript{F}Z^0\right\|^0
\leq\left\|\EuScript{F}Z^0\right\|^0
=\left\|Z^0\right\|^0\leq\left\|Z^0\right\|^1,\  t\in[0,T].
\end{equation}
 Due to \eqref{ttne2} and \eqref{ttne3}, we  obtain
\begin{align}
\Bigg|\int_0^\infty K_{y_1}(0,x)&\Bigg(\frac{e^{-\frac{x^2}{4t}}}{\sqrt{4\pi t}}* Z^0(x)\Bigg)dx\Bigg|
\leq\frac{1}{\sqrt 2}\left\|K_{y_1}(0,\cdot)\right\|_{L^2(0,+\infty)}\Bigg\|\frac{e^{-\frac{(\cdot)^2}{4t}}}{\sqrt{4\pi t}}*Z^0\Bigg\|^0\nonumber\\
&\leq\frac{1}{\sqrt2}\sqrt{2\sigma_0(0)\left(R_0+M_1^2R\right)}\left\|Z^0\right\|^1
\nonumber\\
&=\sqrt{\sigma_0(0)\left(R_0+M_1^2R\right)}\left\|Z^0\right\|^1,\quad t\in[0,T]. \label{ttne4}
\end{align}
 With regard to \eqref{ttne1}, \eqref{ttne4}, and \eqref{est3ne}, we get
\begin{align}
|f(t)|&\leq\sqrt[4]{(\rho k)(0)}\|u^{\rho k\gamma}\|_{L^\infty(0,T)}+\frac{1}{2\sqrt{2}}\sigma_0(0)\left\|Z^0\right\|^1
\nonumber\\
&\kern2.8ex+\sqrt{\sigma_0(0)\left(R_0+M_1^2R\right)}\left\|Z^0\right\|^1
\nonumber\\
&\leq\sqrt[4]{(\rho k)(0)}\|u^{\rho k\gamma}\|_{L^\infty(0,T)}
\nonumber\\
&\kern2.8ex+E_2\left(\frac{\sigma_0(0)}{2\sqrt{2}}+\sqrt{\sigma_0(0)\left(R_0+M_1^2R\right)}\right)\fnl W^0\fnr^1,\quad t\in[0,T].\label{estfne}
\end{align}
Taking into account \eqref{estk1ne}, we obtain
\begin{align}
|P(t)|&\leq\frac{1}{\sqrt{\pi t}}\left(\frac{1}{4}\int_0^\infty\left|r\left(\frac{x}{2}\right)\right|dx+M_1\int_0^\infty\left|\sigma_0\left(\frac{x}{2}\right)\right|dx+\frac{1}{2}\sigma_0(0)\right)\nonumber\\
&=\frac{\sigma_0(0)+2M_1R}{\sqrt{\pi t}},\quad t\in[0,T].\label{estpne}
\end{align}
Using \eqref{estfne} and \eqref{estpne} and applying Lemma \ref{exist}, we conclude that there exists a unique solution to equation \eqref{eqoper} (and, consequently, \eqref{contrneu}). Moreover, using \eqref{estvvv}, we have
\begin{align*}
&\|u^{110}\|_{L^\infty(0,T)}\leq\left(1+2\left(\sigma_0(0)+2M_1R\right)\sqrt{\frac{T}{\pi}}\right)e^{\left(\sigma_0(0)+2M_1R\right)^2T}\\
&\times\left(\sqrt[4]{(\rho k)(0)}\|u^{\rho k\gamma}\|_{L^\infty(0,T)}+E_2\left(\frac{\sigma_0(0)}{2\sqrt{2}}+\sqrt{\sigma_0(0)\left(R_0+M_1^2R\right)}\right)\fnl W^0\fnr^1\right).
\end{align*}
The theorem is proved.
\end{proof}

Due to Theorems \ref{thtransf1ne} and \ref{thtransf2ne},
the operator $\widehat{\mathbb T}$  not only is a continuous one-to-one mapping between the spaces $\widehat H^s$ and $\widehat{\mathbb H}^s$ (see Theorem \ref{th-trans-neu}) but also is one-to-one mapping between the set of the solutions to \eqref{eq}--\eqref{ic}  with constant coefficients ($\rho=k=1$, $\gamma=0$) where $u = u^{110}\in L^\infty(0,T)$ and the set of the solutions to this problem with variable coefficients $\rho, k, \gamma$ where $u=u^{\rho k\gamma}\in L^\infty(0,T)$, where $u^{110}$ and $u^{\rho k\gamma}$ are different generally speaking.

In \cite[Section 7]{FKh2},  piecewise constant controls $u_{N,l}^{110}$, $N,l\in\mathbb{N}$, solving the  approximate controllability problem for system \eqref{z1ne}, \eqref{z2ne}, have been constructed.
Moreover, the solution to this system with the controls $u_{N,l}^{110}$ has been  obtained:
\begin{align*}
Z_{N,l}(\xi,t)=\frac{e^{-\frac{\xi^2}{4t}}}{2\sqrt{\pi t}}&\ast Z^0(\xi)
\\
&-\sqrt{\frac{2}{\pi}}\int\limits_0^t e^{-\frac{\xi^2}{4\tau}}\frac{u_{N,l}^{110}(t-\tau)}{\sqrt{2\tau}}d\tau,\quad N,l\in\mathbb{N},\ \xi\in\mathbb{R},\ t\in[0,T].
\end{align*}
In addition, it has been proved that
$$
\| Z^T-Z_{N,l}(\cdot,T)\|^1\rightarrow 0 \quad \text{as}\ N\rightarrow\infty\ \text{and}\ l\rightarrow\infty.
$$

Therefore, according to  Theorem \ref{thtransf1ne}, the controls 
\begin{align*}
u_{N,l}^{\rho k\gamma}(t)&=\frac{1}{\sqrt[4]{(\rho k)(0)}}\left(u_{N,l}^{110}(t)+\int\limits_0^\infty K_{y_1}(0,\xi)Z_{N,l}(\xi,t)d\xi
\right.\\
&\kern15ex\left.-\frac{1}{2}Z_{N,l}(0^+,t)\int\limits_0^\infty r(\xi)d\xi\right),\quad t\in[0,T],\ N\in\mathbb{N},\ l\in\mathbb{N},
\end{align*}
 solve the  approximate controllability problem for system \eqref{eq1ne}, \eqref{ic1ne} with $u=u_{N,l}^{\rho k\gamma}$. 
In addition, $u_{N,l}^{\rho k\gamma}\in L^\infty(0,T)$ due to Theorem \ref{thtransf1ne}. Moreover,
$W_{N,l}(\cdot,t)=\widehat{\mathbb{T}}Z_{N,l}(\cdot,t)$, $t\in[0,T]$, and 
$$
\fnl W^T-W_{N,l}(\cdot,T)\fnr^1\rightarrow 0 \quad \text{as}\ N\rightarrow\infty\ \text{and}\ l\rightarrow\infty.
$$

\section{Examples}
\label{exxx}

\begin{example}
\label{exxx2}
Consider  system \eqref{eq}--\eqref{ic} with 
\begin{align*}
k(x)&=\frac{(1+2|x|)\cosh x}{3},\quad \rho(x)=\frac{12\cosh x}{1+2|x|},
\\ 
\gamma(x)&=\frac{(1+2|x|)\tanh |x|}{36}+\frac{(1+2|x|)^2}{144}\left(1+\frac{1}{\cosh^2 x}\right)-\frac{1}{4(1+2|x|)^3},\quad x\in\mathbb{R}.
\end{align*}
We have
\begin{align*}
Q_2(\rho,k)&=\frac{(1+2|x|)\tanh |x|}{36}+\frac{(1+2|x|)^2}{144}\left(1+\frac{1}{\cosh^2 x}\right),\\ q(x)&=Q_2(\rho,k)-\gamma(x)=\frac{1}{4(1+2|x|)^3},\quad x\in\mathbb{R}.
\end{align*}
Due to \eqref{sigma}, we get 
$$
\sigma(x)=\operatorname{\mathrm{sgn}} x\ln\left(1+2|x|\right)^3,\ x\in\mathbb{R},\quad\text{and}\quad \sigma^{-1}(\lambda)=\frac{1}{2}\operatorname{\mathrm{sgn}} \lambda \left(e^{\frac{|\lambda|}{3}}-1\right),\ \lambda\in\mathbb{R}.
$$

Let us  consider system \eqref{z1ne}, \eqref{z2ne} with  
$Z^0(x)=e^{-\frac{|x|}{2}}$ and with the steering condition $Z^T(x)=e^{-\frac{2|x|-T}{4}}$, $x\in\mathbb R$.
Evidently, 
$$
Z(x,t)=e^{-\frac{2|x|-t}{4}},\quad x\in\mathbb R,\ t\in[0,T],
$$
is the unique solution to this system and the state $Z^0$ is controllable to the state $Z^T$ with respect to system \eqref{z1ne}, \eqref{z2ne} in the time $T$ with the control 
$$
u^{110}(t)=-\frac{1}{2}e^\frac{t}{4},\quad t\in[0,T].
$$

Now consider  system \eqref{eq1ne}, \eqref{ic1ne} with the given  $q$.
According to Theorem \ref{corol2ne} \ref{daa1ne},
the state  $W^0=\widehat{\mathbb{T}}Z^0$ is controllable to the state $W^T=\widehat{\mathbb{T}}Z^T$ with respect to system \eqref{eq1ne}, \eqref{ic1ne} in the time $T$. Moreover, due to Theorem \ref{thtransf1ne}, a control $u^{\rho k\gamma}$ solving controllability problem for system \eqref{eq1ne}, \eqref{ic1ne} is defined by \eqref{contr1ne}.

Let us find $W^0$, $W^T$ and $u^{\rho k\gamma}$ explicitly.
Due to \eqref{r}, $r(\lambda)=q\circ\sigma^{-1}=\frac{1}{4}e^{-\lambda}$, $\lambda>0$. The kernel of the transformation operator $\widehat{\textbf{T}}_r$ has been found in \cite{F2,Kh2} for this $r$. We have
\begin{equation}
\label{kkker}
K(y_1,y_2)=\frac{e^{-\frac{y_1+y_2}{2}}}{4}\frac{I_1\left(\sqrt{e^{-\frac{y_1}{2}}\left(e^{-\frac{y_1}{2}}-e^{-\frac{y_2}{2}}\right)}\right)}{\sqrt{e^{-\frac{y_1}{2}}\left(e^{-\frac{y_1}{2}}-e^{-\frac{y_2}{2}}\right)}},\quad y_2>y_1>0,
\end{equation}
where $I_n$ is the modified Bessel function, $n=\overline{0,\infty}$.

It is well-known (see, e.g., \cite[9.6.28]{HB}) that \begin{equation}
\label{abr-st}
\big(y^{-n}I_n(y)\big)'= y^{-n}I_{n+1}(y)\  \text{and}\  \big(y^{n}I_n(y)\big)'= y^{n}I_{n-1}(y),\quad y>0,\ n\in\mathbb N.
\end{equation}
Due to the first of these formulae with $n=1$, we obtain
\begin{align*}
&K_{y_1}(y_1,y_2)=-\frac{1}{8}e^{-\frac{y_1+y_2}{2}}\frac{I_1\left(\sqrt{e^{-\frac{y_1}{2}}\left(e^{-\frac{y_1}{2}}-e^{-\frac{y_2}{2}}\right)}\right)}{\sqrt{e^{-\frac{y_1}{2}}\left(e^{-\frac{y_1}{2}}-e^{-\frac{y_2}{2}}\right)}}\\
&-\frac{1}{16}e^{-\frac{y_1+y_2}{2}}\frac{I_2\left(\sqrt{e^{-\frac{y_1}{2}}\left(e^{-\frac{y_1}{2}}-e^{-\frac{y_2}{2}}\right)}\right)}{\sqrt{e^{-\frac{y_1}{2}}\left(e^{-\frac{y_1}{2}}-e^{-\frac{y_2}{2}}\right)}}\frac{e^{-\frac{y_1}{2}}\left(2e^{-\frac{y_1}{2}}-e^{-\frac{y_2}{2}}\right)}{\sqrt{e^{-\frac{y_1}{2}}\left(e^{-\frac{y_1}{2}}-e^{-\frac{y_2}{2}}\right)}},\quad y_2>y_1>0.
\end{align*}
Hence, 
\begin{align*}
K_{y_1}(0,x)&=-\frac{e^{-\frac{x}{2}}}{8}\frac{I_1\left(\sqrt{1-e^{-\frac{x}{2}}}\right)}{\sqrt{1-e^{-\frac{x}{2}}}}
\\
&\kern2.8ex-\frac{e^{-\frac{x}{2}}}{16}I_2\left(\sqrt{1-e^{-\frac{x}{2}}}\right)\left(1+\frac{1}{1-e^{-\frac{x}{2}}}\right),\quad x>0.
\end{align*}
Taking into account \eqref{contr1ne}, we get
\begin{align*}
u^{\rho k\gamma}(t)&=-\frac{e^\frac{t}{4}}{8\sqrt{2}}\left(5+\int_0^\infty e^{-x}\frac{I_1\left(\sqrt{1-e^{-\frac{x}{2}}}\right)}{\sqrt{1-e^{-\frac{x}{2}}}}\,dx\right.\nonumber\\
&\kern2.6ex\left.+\frac{1}{2}\int_0^\infty e^{-x}I_2\left(\sqrt{1-e^{-\frac{x}{2}}}\right)\left(1+\frac{1}{1-e^{-\frac{x}{2}}}\right)dx\vphantom{\frac{I_1\left(\sqrt{1-e^{-\frac{x}{2}}}\right)}{\sqrt{1-e^{-\frac{x}{2}}}}}\right),\quad t\in[0,T].
\end{align*}
Substituting $y$ for $\sqrt{1-e^{-\frac{x}{2}}}$ and then integrating the first integral by parts, we get
\begin{align*}
%\label{contfor}
u^{\rho k\gamma}(t)&=-\frac{e^{\frac t4}}{8\sqrt 2}
\left(5+4\int_0^1\big(1-y^2\big)I_1(y)\, dy
+2\int_0^1\left(\frac1y-y^3\right)I_2(y)\, dy \right)
\nonumber\\
&= -\frac{e^{\frac t4}}{8\sqrt 2} \left(1+8\int_0^1 yI_0(y)\, dy +2\int_0^1 \frac{I_2(y)}{y}\, dy -2 \int_0^1 y^3I_2(y)\, dy\right).
\end{align*}
With regard to \eqref{abr-st}, we obtain 
\begin{align}
\label{contrfin}
u^{\rho k\gamma}(t)
&=-\frac{e^{\frac t4}}{8\sqrt 2}\left(1+8I_1(1)+2I_1(1)-1-2I_3(1)\right)
\notag \\
&=\frac{e^\frac{t}{4}}{4\sqrt{2}} \left( I_3(1)-5I_1(1)\right),\quad t\in[0,T].
\end{align}
According to the definition of the operator $\widehat{\mathbb{T}}$, we  have
\begin{equation*}
W(x,t)=\left(\textbf{S}\widehat{\textbf{T}}_rZ\right)(x,t)=\frac{\left(\widehat{\textbf{T}}_rZ(\cdot,t)\right)\left(\operatorname{\mathrm{sgn}} x\ln\left(1+2|x|\right)^3\right)}{\sqrt{2\cosh x}},\  x\in\mathbb{R},\ t\in[0,T].
\end{equation*}
Due to the definition of the operator $\widehat{\textbf{T}}_r$, we obtain
\begin{align*}
\left(\widehat{\textbf{T}}_rZ\right)&(\lambda,t)=e^{-\frac{2|\lambda|-t}{4}}
\\
&+\int_{|\lambda|}^\infty \frac{e^{-\frac{|\lambda|+x}{2}}}{4}\frac{I_1\left(\sqrt{e^{-\frac{|\lambda|}{2}}\left(e^{-\frac{|\lambda|}{2}}-e^{-\frac{x}{2}}\right)}\right)}{\sqrt{e^{-\frac{|\lambda|}{2}}\left(e^{-\frac{|\lambda|}{2}}-e^{-\frac{x}{2}}\right)}}e^{-\frac{2x-t}{4}}dx,\ \lambda\in\mathbb{R}, \ t\in[0,T].
\end{align*}
Replacing $\sqrt{e^{-\frac{|\lambda|}{2}}\left(e^{-\frac{|\lambda|}{2}}-e^{-\frac{x}{2}}\right)}$ by $y$ in the integral, we get
\begin{align}
\left(\widehat{\textbf{T}}_rZ\right)&(\lambda,t)=e^{-\frac{2|\lambda|-t}{4}}+ e^\frac{2|\lambda|+t}{4}\int_0^{e^{-\frac{|\lambda|}{2}}}(e^{-|\lambda|}-y^2)I_1(y)dy=e^{-\frac{2|\lambda|-t}{4}}\nonumber\\
&+ e^\frac{t}{4}\left(2I_1\left(e^{-\frac{|\lambda|}{2}}\right)-e^{-\frac{|\lambda|}{2}}\right)=2e^\frac{t}{4} I_1\left(e^{-\frac{|\lambda|}{2}}\right),\quad \lambda\in\mathbb{R},\ t\in[0,T].\label{ttilde}
\end{align}
Thus,
\begin{align*}
W(x,t)&=e^\frac{t}{4}\sqrt{\frac{2}{\cosh x}}I_1\left(e^{-\frac{1}{2}\ln\left(1+2|x|\right)^3}\right)
\\
&=e^\frac{t}{4}\sqrt{\frac{2}{\cosh x}}I_1\left(\frac{1}{\left(1+2|x|\right)^{3/2}}\right),\ x\in\mathbb{R},\ t\in[0,T].
\end{align*}
Hence,
\begin{align}
W^0(x)&=\sqrt{\frac{2}{\cosh x}}I_1\left(\frac{1}{\left(1+2|x|\right)^{3/2}}\right),\quad x\in\mathbb{R},\label{initne}\\
W^T(x)&=e^\frac{T}{4}\sqrt{\frac{2}{\cosh x}}I_1\left(\frac{1}{\left(1+2|x|\right)^{3/2}}\right),\quad x\in\mathbb{R}.\label{enddne}
\end{align}
Thus, the initial state $W^0$ defined by \eqref{initne} is controllable to the steering state $W^T$ defined by \eqref{enddne} with respect to system \eqref{eq1ne}, \eqref{ic1ne} in the time $T$ by the control \eqref{contrfin}.
\end{example}

%%%%%%%%%%%%%%%

\begin{figure}[!h]
	\begin{center}
		\begin{subfigure}[t]{0.49\linewidth}
			\centering \includegraphics[height=45mm]{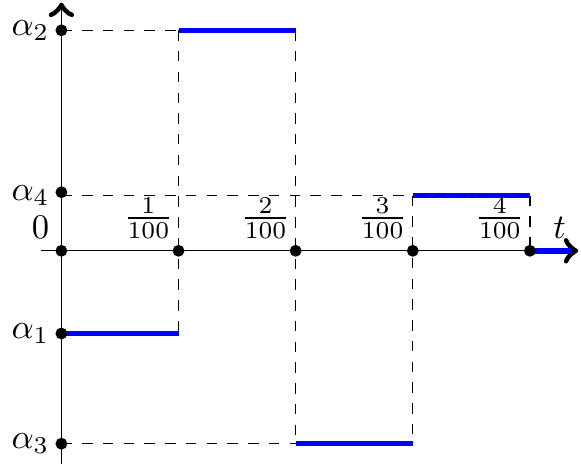}\\
			\centering \caption{\parbox[t]{0.88\textwidth}{$N=3$, $l=100$,\\
					$\mathrm \alpha_1\approx -119704.546455$,\\
					$\mathrm \alpha_2\approx 318558.179365$,\\
					$\mathrm \alpha_3\approx -282251.95269$,\\
					$\mathrm \alpha_4\approx 83317.88255$.}}
			\label{fig:contr1}
		\end{subfigure}
		\begin{subfigure}[t]{0.49\linewidth}
			\centering \includegraphics[height=46mm]{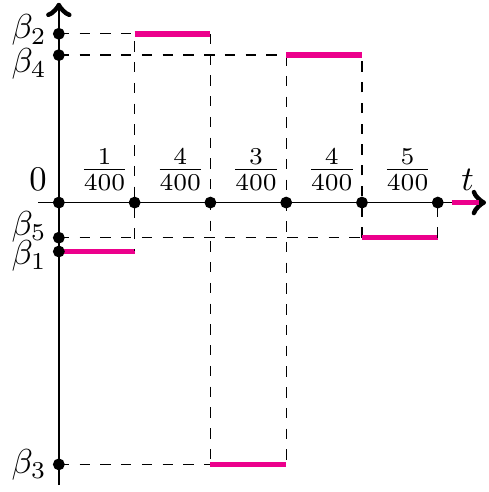}
			\centering \caption{\parbox[t]{0.88\textwidth}{%
					$N=4$, $l=400$, \\
					$\mathrm \beta_1\approx -183378505.929335$,\\
					$\mathrm \beta_2\approx 701420689.4293751$, \\
					$\mathrm \beta_3\approx -1006324503.657385$, \\
					$\mathrm \beta_4\approx 641835320.740755$,\\
					$\mathrm \beta_5\approx -153553322.43498$.}}
			\label{fig:contr2}
		\end{subfigure}
	\end{center}
	\centering \caption{The controls $u_{N,l}^{\rho k\gamma}$.}
	\label{fig:contr}
\end{figure}

\begin{figure}[!h]
	\begin{center}
		\begin{subfigure}[t]{0.49\linewidth}
			\centering \includegraphics[height=45mm]{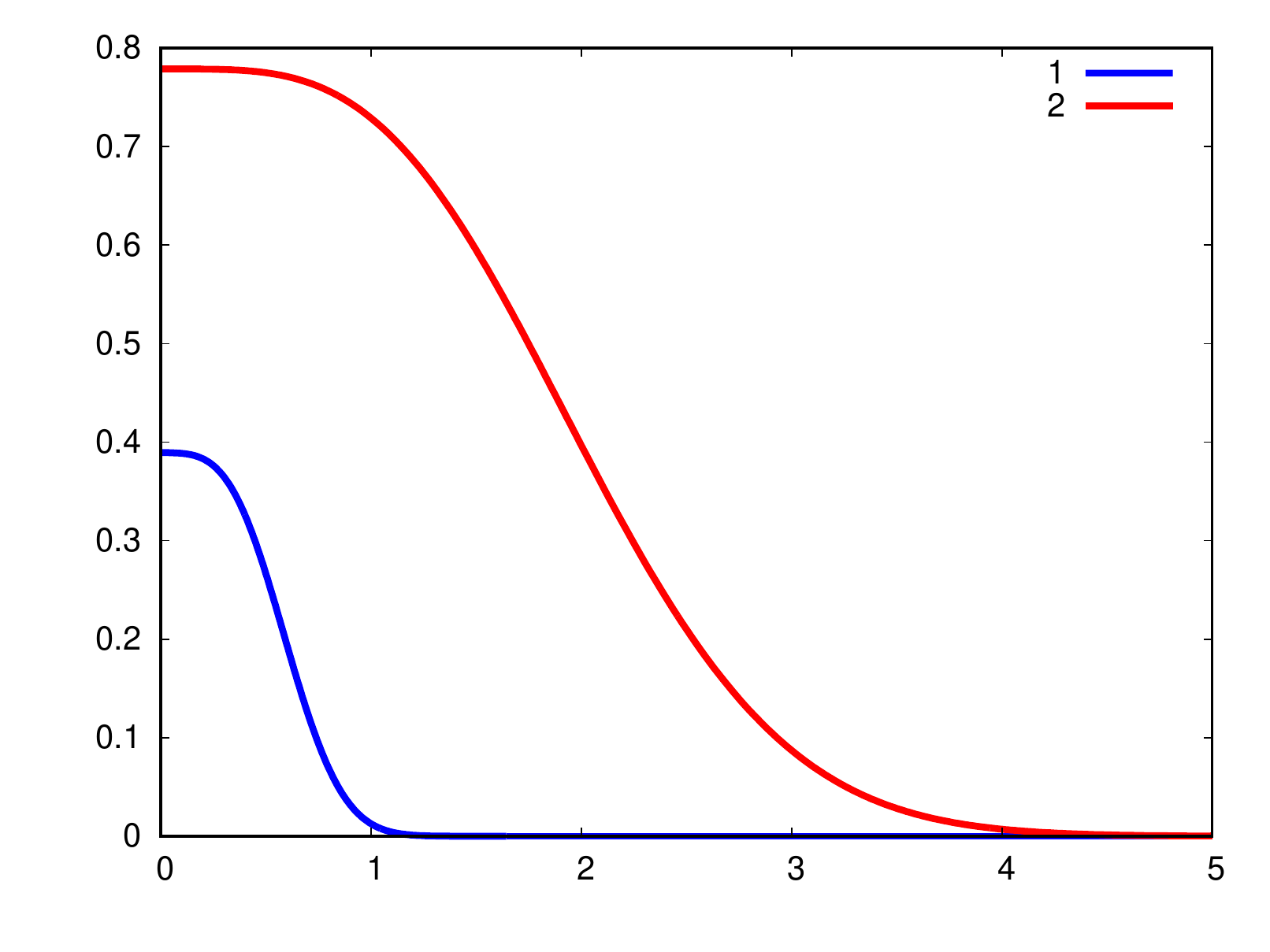}\\
			\centering \caption{\parbox[t]{0.88\textwidth}{\zcx{1}~The given target state $W^T$;
					\zcx{2}~The function $Z^T=\widehat{\mathbb{T}}^{-1}W^T$.\\  \mbox{}\\ \mbox{}\\ \mbox{}}}
			\label{fig:appr11}
		\end{subfigure}
		\begin{subfigure}[t]{0.49\linewidth}
			\centering \includegraphics[height=46mm]{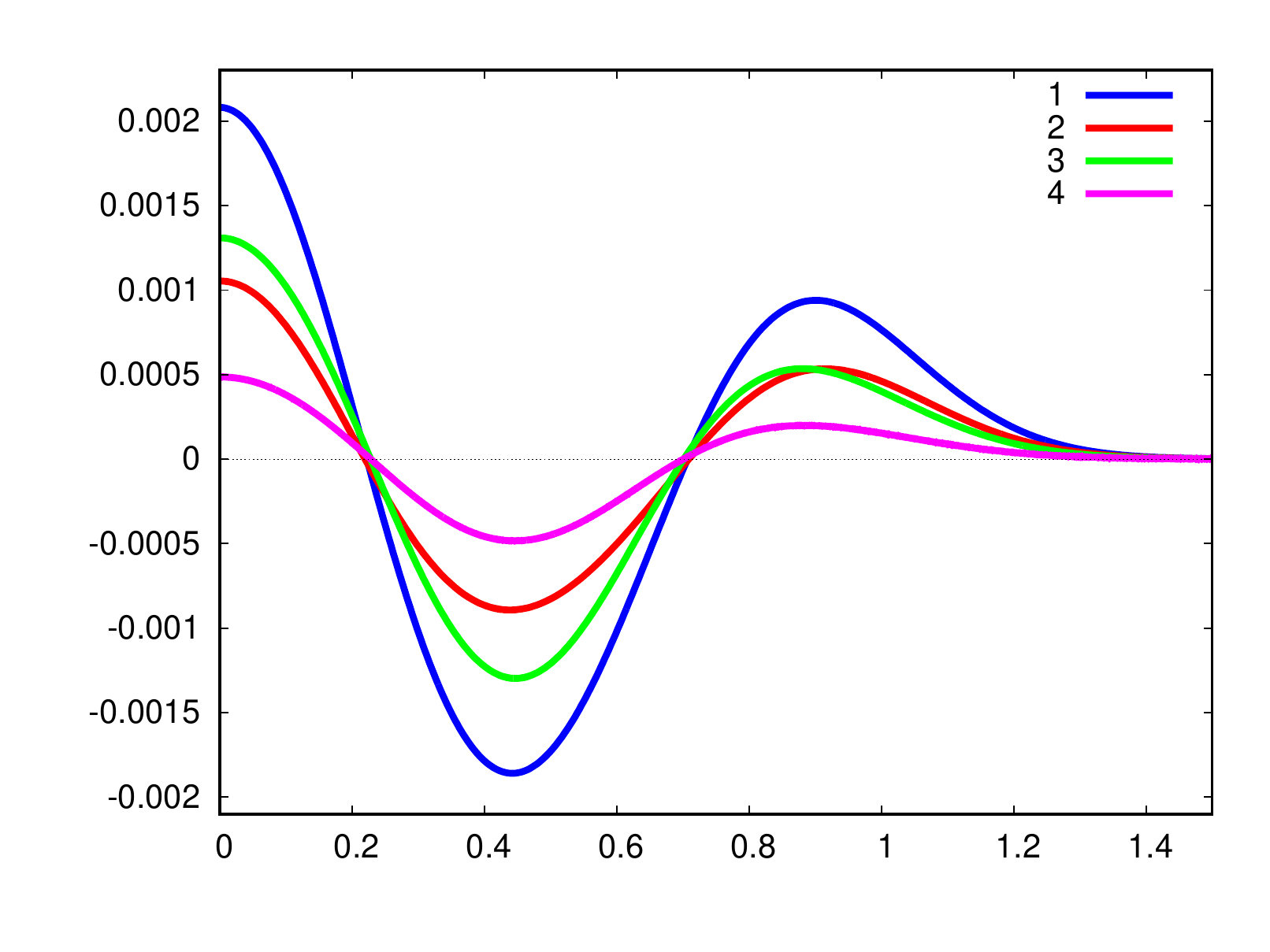}
			\centering \caption{\parbox[t]{0.88\textwidth}{%
					The difference $W^T-W_{N,l}^T$ in the cases: \zcx{1}~$N=3$, $l=100$; \zcx{2}~$N=3$, $l=200$; \zcx{3}~$N=4$, $l=150$; \zcx{4}~$N=4$, $l=400$.}}
			\label{fig:appr12}
		\end{subfigure}
	\end{center}
	\centering \caption{The influence of the controls $u=u_{N,l}^{\rho k\gamma}$
		on the end state $W_{N,l}^T$ of the solution to \eqref{eq1ne}, \eqref{ic1ne}.}
	\label{fig:appr2}
\end{figure}

\begin{example}
\label{exxx3}
Let 
$$
k(x)=\frac{4+x^2}{3+|x|},\quad \rho(x)=(4+x^2)(3+|x|),\quad \gamma(x)=\frac{12-|x|^3}{(3+|x|)^3(4+x^2)^2},\quad x\in\mathbb R.
$$
Consider approximate controllability problem   for system \eqref{eq}--\eqref{ic}, \eqref{ec}, where $T=1/2$, $w^0=0$, $u=u^{\rho k\gamma}$, and
\begin{equation*}
w^T(x)=\frac{1}{\sqrt{4+x^2}}\cosh\frac{x(|x|+6)}{2\sqrt{2T}}e^{-\frac{x^2\left(|x|+6\right)^2}{16T}-\frac{1}{4}},\quad x\in\mathbb R.
\end{equation*}
It is easy to see that 
$$
Q_2(\rho,k)=\frac{12-|x|^3}{(3+|x|)^3(4+x^2)^2},\quad x\in\mathbb{R}.
$$ 
Therefore, $q(x)=Q_2(\rho,k)-\gamma(x)=0$ on $\mathbb{R}$.
We obtain 
$$
\sigma(x)=\frac{1}{2}x\left(|x|+6\right),\quad x\in\mathbb{R},\quad\text{and}\quad \sigma^{-1}(\lambda)=\operatorname{\mathrm{sgn}} \lambda \left(\sqrt{2|\lambda|+9}-3\right),\quad \lambda\in\mathbb{R}.
$$
We have $W^0=w^0$ and $W^T=w^T$ on $\mathbb{R}$. Consider control system \eqref{eq1ne}, \eqref{ic1ne} with $q=0$, $W^0=0$, and with the steering condition 
$$
W(x,T)=W^T(x)=\frac{1}{\sqrt{4+x^2}}\cosh\frac{x(|x|+6)}{2\sqrt{2T}}e^{-\frac{x^2\left(|x|+6\right)^2}{16T}-\frac{1}{4}},\quad x\in\mathbb{R},\  T=1/2.
$$ 
Let us investigate whether the state $W^0$ is approximately controllable to a target state $W^T$ with respect to system \eqref{eq1ne}, \eqref{ic1ne} in the time $T=1/2$.
  
According to \eqref{r}, $r=0$ on $\mathbb{R}$. Hence, $\widehat{\textbf{T}}_r=\mathrm{Id}$, and the transformation operator $\widehat{\mathbb{T}}$ takes the form
$\widehat{\mathbb{T}}=\textbf{S}$.
Denote 
$Z(\cdot,t)=\left(\widehat{\mathbb{T}}^{-1}W\right)(\cdot,t)=\left(\textbf{S}^{-1}W\right)(\cdot,t)$, $t\in[0,T]$, $Z^0=\widehat{\mathbb{T}}^{-1}W^0=\textbf{S}^{-1}W^0$, $Z^T=\widehat{\mathbb{T}}^{-1}W^T=\textbf{S}^{-1}W^T$. 

Due to Theorem \ref{thtransf2ne}, $Z$ is the solution to system \eqref{z1ne}, \eqref{z2ne} with 
\begin{align*}
u=u^{110}=\sqrt[4]{(\rho k)(0)}u^{\rho k\gamma}=2u^{\rho k\gamma},\quad Z^0=0,
\end{align*}
and with the steering condition
\begin{align*}
Z(\xi,T)=Z^T(\xi)=\cosh\frac{\xi}{\sqrt{2T}} e^{-\frac{\xi^2}{4T}-\frac{1}{4}},\quad \xi\in\mathbb{R},\ T=1/2.
\end{align*}
Controllability problems  for this system have been considered in Example 4 in \cite{FKh2}.
Controls solving the approximate controllability problem
for system \eqref{z1ne}, \eqref{z2ne} have been found in the form
$$
u_{N,l}^{110}=\sum_{p=0}^N U_{p,l}^N,\quad N\in\mathbb{N},
$$
where $U_{p,l}^N\in \mathbb{R}$ is a constant, $p=\overline{0,N}$, $l$ depends on $N$, $N\in\mathbb N$.
The end states $Z_{N,l}^T$ such that 
$$
\forall \varepsilon>0 \ \exists N\in\mathbb N\ \exists l\in\mathbb N\quad \left\| Z^T-Z_{N,l}^T\right\|^1\leq\varepsilon
$$ 
have been found in the form
\begin{equation*}
Z_{N,l}^T(\xi)=-\sqrt{\frac2\pi}\int_0^T e^{-\frac{\xi^2}{4\tau}}\frac{u_{N,l}^{110}(T-\tau)}{\sqrt{2\tau}}d\tau,\quad\xi\in\mathbb{R},\ T=1/2.
\end{equation*}
Applying Theorem \ref{thtransf1ne}, we conclude that 
controls $u_{N,l}^{\rho k\gamma}=\frac{1}{2}u_{N,l}^{110}$ solve the approximate controllability problem for given system \eqref{eq1ne}, \eqref{ic1ne}. Moreover,
\begin{equation*}
W_{N,l}^T(x)=\frac{1}{\sqrt{4+x^2}}Z_{N,l}^T\left(\frac{1}{2}x\left(|x|+6\right)\right),\quad x\in\mathbb{R},\ T=1/2,
\end{equation*}
and for $\varepsilon$, $N$, and $l$ mentioned above, we get
\begin{equation*}
\fnl W^T-W_{N,l}^T\fnr^1\leq E_0\varepsilon,
\end{equation*}
where $E_0$ is the constant from estimate \eqref{est1ne}.
The graphs of $u_{N,l}^{\rho k\gamma}$ and $W_{N,l}^T$  see in Figs. \ref{fig:contr}, \ref{fig:appr2}.
\end{example}


\section{Acknowledgements}

The authors are partially supported by ``Pauli Ukraine Project'', funded in the WPI Thematic Program ``Quantum equations and experiments'' (2021/2022). The second author is partially supported by  the National Academy of Sciences of Ukraine (Grant No. 0122U111111).

%\clearpage

%%%%%%%%%%%%%%%%%%%%%%%%%%%%%%

\EndPaper%%%%  do not remove this command

%%%%%%%%%%%%%%%%%%%%%%%%%%%%%%%%%%%%%%%%%%%%%%%%%%%%%%%

\end{document}